\documentclass[12pt]{amsart}
\usepackage{amsmath,amssymb,amscd}
\usepackage{epsf}
\newtheorem{prop}{Proposition}%[section]
\newtheorem{thm}[prop]{Theorem}

\newtheorem{lem}[prop]{Lemma}
\def\demo#1{\medskip\noindent{\bf#1\enspace}}
\def\enddemo{\medskip}
\def\QED{\hbox{\hskip 6pt\vrule width6pt height7pt depth1pt \hskip1pt}}

        \def\zed{{\mathbb Z}}
        \def\red{{\mathbb R}}
        
        \def\rat{{\mathbb Q}}
        
\begin{document}

\title{Tiling spaces are Cantor set fiber bundles}

\author{Lorenzo Sadun and R. F. Williams}

%\dedicatory{}

\thanks{Mathematics subject classification(1985 revision). Primary
  58F13, 54F15.  Key words:  tiling, tiling space, branched surface,
  Cantor fiber bundle.} 
\thanks{R.F.W gratefully acknowledges helpful conversations with
  Marcy Barge and hospitality of the Mathematics Department of Montana
  State University. This work is partially supported by the National
Science Foundation and the Texas Advanced Research Program}

\begin{abstract}We prove that fairly general spaces of tilings of
$\red^d$ are fiber bundles over the torus $T^d$, with totally
disconnected fiber.  This was conjectured (in a weaker form) in [W3],
and proved in certain cases.  In fact, we 
show that each such space is homeomorphic to the  
$d$-fold suspension of a $\zed^d$ subshift (or equivalently, a tiling
space whose tiles are marked unit $d$-cubes). The only restrictions on
our tiling spaces are that 1) the tiles are assumed to be polygons
(polyhedra if $d>2$) that meet full-edge to full-edge (or full-face
to full-face), 2) only a finite number of tile types are allowed, and
3) each tile type appears in only a finite number of orientations.
The proof is constructive, and we illustrate it by constructing a 
``square'' version of the Penrose tiling system.

\end{abstract}

\maketitle

\markboth{Sadun-Williams}{Tiling spaces are Cantor set bundles}

\section{Introduction and results}
Anderson and Putnam [AP] have found a relation between tiling theory
(specifically, substitution tiling theory) and a theory of attractors
studied in the dynamics of diffeomorphisms of manifolds with
hyperbolic structure.  They show that {\it substitution} tiling spaces
are a special case of expanding attractors, a concept introduced in
[W2] to study the dynamics of diffeomorphisms.  It is well known that
an expanding attractor is locally the topological product of a Cantor
set and a disk of the appropriate dimension.  But what is it
globally?  Can the local neighborhoods be stitched together to form a
bundle over a manifold whose fiber is a Cantor set? In 1973, Williams
[W2] conjectured that they could. This conjecture has long been known
to be true in dimension 1, ([BJV], [W1] etc).  However, it fails in
higher dimensions; Farrell and Jones [FJ] constructed a
counterexample in dimension 2.

In this paper, however, we show (Theorem 1) that the conjecture is
true for tiling spaces in all dimensions, even without the
substitution assumption. Tiling spaces in general are {\it flat},
suggesting that the appropriate base spaces should be tori. 
(In fact, the torus is the only candidate for any expanding attractor
since we know that the area must grow as a polynomial of the 'radius'.)
In addition, we show (Theorem 2) that tiling spaces are
suspensions of $\zed^d$ subshifts, which are manifestly bundles over tori.

The eventual goal is to classify tiling spaces up to homeomorphism,
and the existence of a bundle structure is an important first step.
Although Cantor set fiber bundles are harder to deal with than, say,
vector bundles, topological invariants have been constructed for bundles
over $S^1$ ([BF,BJV,F,S]).  In higher dimensions, 
the bundles are generically less
``flabby'' and thus, hopefully, allow stronger invariants than 
those of Parry-Sullivan [PS] and Bowen-Franks [BF] type.
In addition, Theorem 2 suggests that much of the theory of subshifts [B,LM,Wa]
can be brought into play.

The proofs given below are easy, and are too abstract and too general
to be of direct use in classification. For note that given a bundle
structure $p:X \to T^d$ and an $n$-to-1 self covering $f: T^d \to T^d$,
the composition $f \circ p: X \to T^d$ is another bundle structure.
The latter is ``coarser'' in that each of its fibers contains $n$ fibers
of the former; is there a finest possible bundle structure?  The current
construction does not address this problem.  However, see [W3]
for the case of the Penrose tiling.  The 40 tiles of the Penrose
tiling cover the torus 10 times in [W3], but 232 times in the current
version.

The basic ingredients of this paper are tiling systems $P$ of
$\red ^d,$ and the corresponding tiling spaces $X(P).$  These
spaces are assumed to satisfy the following hypotheses:

\begin{enumerate}
\item The tiles are (triangulated) polyhedra that meet full-face to
  full-face. 
\item Only a finite number of tile types (proto-tiles) appear.  In
  this counting, 
tiles that are translations of one another are considered to be the
same type, but tiles that are rotations of one another are considered
to be different.
\item The space $X(P)$ is a closed, nonempty and translation-invariant 
subset of the space of 
all tilings that can be formed from the tiles in $P$.
\end{enumerate}
We will henceforth refer to both the tiling system and the associated
topological space by the same letter $P$.

Our first result is:
\begin{thm} A tiling space that satisfies the above hypotheses is a 
fiber bundle over the torus, with totally disconnected fiber.
\end{thm}

Note that we do not assume that the tilings are quasi periodic, or 
generated by a substitution, or even that they are non-periodic. The
only difference between these cases is the nature of the fiber.  
The fiber for a substitution tiling, or a quasi periodic tiling, will
be a Cantor set, while the fiber for a ($d$-fold) 
periodic tiling system will be a finite collection of points. 

The requirement of polygonal tiles is mostly for convenience.  
A tiling, such as the Penrose chickens, whose edges follow standard 
shapes, can be deformed to a tiling system with polygonal tiles, and therefore
is a fiber bundles over a torus. The requirement that tiles appear
in only a finite number of orientations is more serious.  The techniques
of this paper do not apply to pinwheel-like tiling spaces [R].

Our second result is:
\begin{thm} A tiling space $P$ that satisfies the above hypotheses is
homeomorphic to a tiling space $S$ whose tiles are marked $d$-cubes, or
equivalently to the $d$-fold suspension of a $\zed^d$ subshift. The
space $S$ is defined by local matching rules if and only if $P$ is. 
\end{thm}

Note that this theorem proves the existence of a homeomorphism, not a
topological conjugacy.  The homeomorphism typically does {\it not}
commute with translations, much less with rotations. For more 
information on $\zed^d$ subshifts, see \cite{Wa}.

The proofs proceed as follows.  We call a tiling space rational
(integral) if each edge of each tile is given by a vector with
rational (integral) coordinates.  In Section 2 we show that every
tiling space $P$ can be deformed to a rational tiling space $R$.  This
deformation is a homeomorphism of tiling spaces, but not a topological
conjugacy.  We then show that every rational tiling is a fiber bundle
over the torus.  This proves Theorem 1.

In Section 3 we prove Theorem 2.  We rescale the rational tiling space
$R$ into an integral tiling space, and replace the straight edges with
zig-zags consisting of unit segments in the several coordinate
directions.  The faces then become unions of unit squares, the 3-cells
become unions of unit cubes, and so on.  This gives a ``zig-zag''
system $Z$. The tiles of $Z$ may take on odd shapes, and may even be
disconnected, but are unions of $d$-cubes.  The space $Z$ is
homeomorphic (topologically conjugate, in fact) to the rescaled $R$.
Finally, we consider each constituent $d$-cube of a tile $z$ in $Z$ to
be a tile in a tiling space $S$, with the matching rule that wherever
one such constituent appears, the other constituents of $t$ also
appear nearby. $S$ is a suspension of a subshift, but is also
topologically conjugate to $Z$, and therefore homeomorphic to $P$.

\newpage
\section{Tiling spaces as fiber bundles}

\begin{lem} A tiling space $P$ meeting the above hypotheses is homeomorphic
to a rational tiling space $R$. Furthermore, $R$ has finite type if and
only if $P$ does.
\end{lem}

\demo{Proof.}  For greater clarity, we go through the proof in
dimension 2 and illustrate how each step applies to the
Penrose system. We then indicate how the proof applies, with small changes,
in any dimension.  

Let the tiles of a tiling space $P$ be represented by polygonal disks
$C_i, i=1,...,c,$ in the plane;  let each edge of each tile be given a
fixed orientation.  If tiles $C_i$ and $C_j$ can meet along a common
edge in a tiling, then these edges must be parallel, of the same
length, and have the same orientation.  That is, the displacement
vectors of these edges 
must be equal.  Thus all together, we have a finite set of these
vectors, say $v_1,...,v_n,$ which we will call {\it edge vectors}.

\begin{figure}
\vbox{
%\centerline{\epsfysize=9.5cm \epsfbox{ssft/pentri} }}
\centerline{\epsfysize=9.5cm \epsfbox{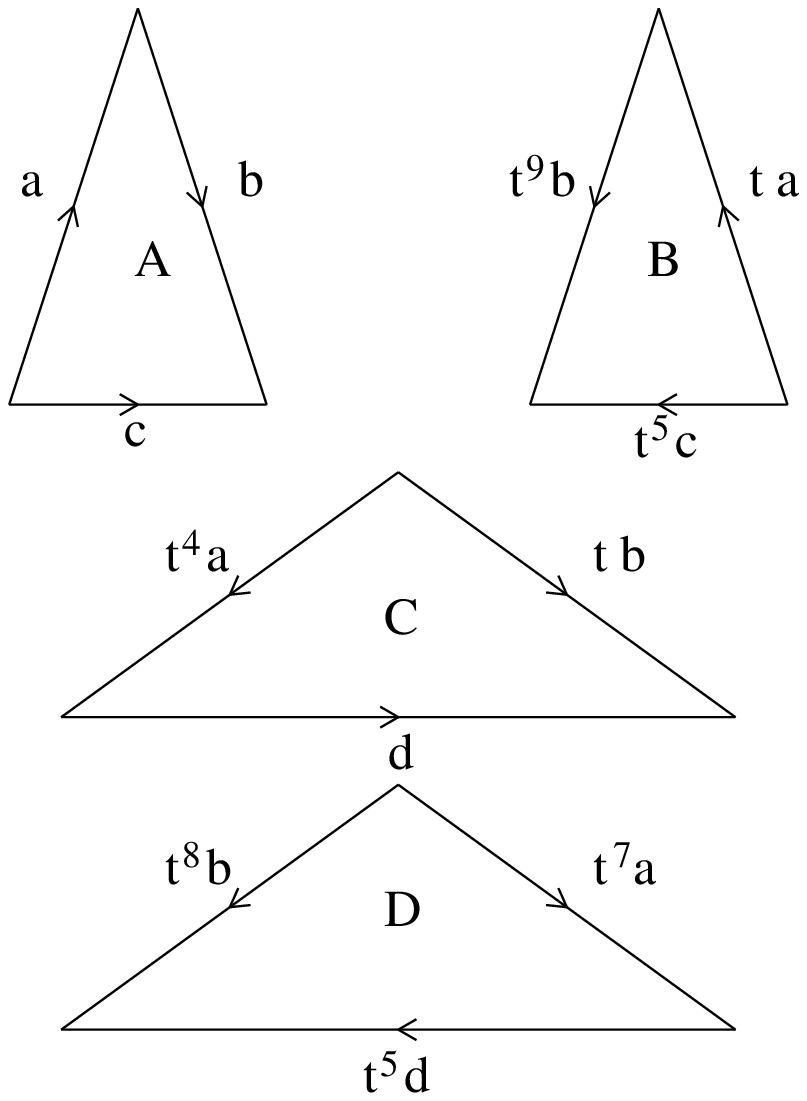} }}
\caption{The tiles of the Penrose system.}
\end{figure}

In the Penrose ``B-tile'' system, there are forty triangular tiles, 
namely those shown in Figure 1 and their rotations by multiples of $2\pi/10$.
We let $t$ denote the rotation by $2\pi/10$, so $t^5A$ means tile $A$
rotated by $\pi$.  Although $A$ is congruent to $B$, they are 
considered separate tiles.  Similarly, $C$ and $D$ are considered distinct.
Because of the identifications, there are 40 edge vectors, not 120.  Five of
the edge vectors are
\begin{eqnarray}
a & = & (2(\tau-1),2\sqrt{\tau+2}) \cr
ta & = &  (-2(\tau-1), 2\sqrt{\tau+2}) \cr
t^2a & = & (-2\tau,2\sqrt{3-\tau}) \cr
t^3a & = & (-4,0)  \cr
t^4a & = & (-2\tau,- 2\sqrt{3-\tau}), 
\end{eqnarray}
where $\tau=(1+\sqrt{5})/2$.  The other vectors are given by
\begin{eqnarray}
t^{5+n}a =&  - t^na & n=0,\ldots 4\cr
t^nb =& t^{n-4} a &n=0, \ldots 9 \cr
t^nc = & (\tau-1)t^{n-2}a& n=0, \ldots 9 \cr
t^n d=& \tau t^{n-2}a & n=0, \ldots 9.
\end{eqnarray}

We wish to construct new tiles $C'_k$, $k=1,\ldots, c$ so that the $n$
(same number!) edge vectors, $v'_1,...v'_n$ of the new tiles all have
rational coordinates.  Furthermore, if an edge $I_i$ of $C_k$ has edge
vector $v_i$, the corresponding edge $I'_i$ of $C'_k$ has edge vector
$v'_i.$ To see that this can be done, note that there is a single
constraint for each tile, to wit the sum of its edge vectors must be
0.  Thus we have a finite system of homogeneous equations:

\begin{equation}
v_{i_1}+ \cdots + v_{i_k}= 0
\label{eq3}
\end{equation}

As the
coefficients in this system are rational (in fact, integers) it
follows by elementary linear algebra, that the solutions with rational
entries are dense in the 
space of all solutions, which is non-empty, as the original vectors
are solutions.  Thus there are rational vectors
$v_1', \cdots, v_n'$ approximating the original vectors, as close as
we wish, which satisfy the constraints.  

The crux here is that we have finitely many loops among the edge
vectors so that the constraints consist of finitely many equations
with integral coefficients.  In higher dimensions the choice of these
loops is a little different.  In fact, since our tiles are
triangulated, the boundary of each 2-dimensional simplex is a loop,
and there are only finitely many of these, which certainly suffice.
Thus in any dimension, the constraint can be taken to consists of a
finite system of equations of the form (3).

%\begin{equation}
%v_{i_1}+ v_{i_2} + v_{i_3}= 0
%\label{eq4}
%\end{equation}

For the Penrose tiling, our system of equations is
\begin{eqnarray}
v'(t^n a) + v'(t^n b) - v'(t^nc) &= 0, & \qquad n=0, \ldots, 9 \cr
v'(t^{n+6}a) + v'(t^{4+n}b) - v'(t^nc) & = 0,  & \qquad n=0, \ldots, 9 \cr
-v'(t^{n+4} a) + v'(t^{n+1} b) - v'(t^n d)  & = 0,  & \qquad n=0, \ldots, 9 \cr
-v'(t^{n+2} a) + v'(t^{n+3} b) - v'(t^n d)  & = 0,  & \qquad n=0, \ldots, 9 
\label{eq5}
\end{eqnarray}
The following is an integer set of solutions:

\smallskip
\vbox{\hfil
\begin{tabular}{|c|c|c|c|c|} \hline
$n$ & $v'(t^n a)$ &  $v'(t^nb)$ & $v'(t^nc)$ & $v'(t^n d)$ \\ \hline
0 & (1,4) & (1,-4) & (2,0) & (6,0) \\
1 & (-1,4) & (3,-2) & (2,2) & (5,4) \\
2 & (-3,2) & (4,0) & (1,2) & (2,6) \\
3 & (-4,0) & (3,2) & (-1,2) & (-2,6) \\
4 & (-3,-2) & (1,4) & (-2,2) & (-5,4) \\
5 & (-1,-4) & (-1,4) & (-2,0) & (-6,0) \\
6 & (1,-4) & (-3,2) & (-2,-2) & (-5,-4) \\
7 & (3,-2) & (-4,0) & (-1,-2) & (-2,-6) \\
8 & (4,0) & (-3,-2) & (1,-2) & (2,-6) \\
9 & (3,2) & (-1,-4) & (2,-2) & (5,-4) \\ \hline
\end{tabular}\hfill}

\smallskip
\noindent Note that by picking integer solutions, we have broken the 
10-fold rotational symmetry of the Penrose system.  This is to be expected,
as one cannot represent $\zed_{10}$ in $GL(2,\rat)$. We have, however, 
preserved the 2-fold rotational symmetry, the reflection symmetry about
the $x$ axis, and the fact that $v'(t^n a) = v'(t^{n+4}b)$. 

Returning to the general case, there is a linear isomorphism from each
edge $I_i$ to the corresponding $I'_i$.  These can be extended to
homeomorphisms from each $C_k$ to $C_k'$.  (If $C_k$ is a triangle
there is a linear extension.  In general, we can always find a
continuous extension.)  We now use these homeomorphisms to convert an
arbitrary tiling by the tiles $\{C_k\}$ into a tiling by the tiles
$\{C'_k\}$.  As we shall see, this procedure is continuous and has a
continuous inverse, and so defines a homeomorphism between the tiling
space $P$ and a rational tiling space $R$.

For each tiling $t$ in the tiling space $P$, we construct a
corresponding tiling $t' \in R$, beginning at the origin.  The origin
in $t$ sits at a point in a closed tile $C_k$; we let the origin in
$t'$ sit at the corresponding point in $C'_k$.  If the origin lies on
an edge $I_i$ that is shared by $C_j$ and $C_k$, then we may start
with either $C_j$ or $C_k$. There is no ambiguity in the location of
the origin, since the maps $C_j \to C_j'$ and $C_k \to C_k'$ are both
extensions of the linear map from $I_i$ to $I_i'$.

This procedure determines the position of a single tile $C_k'$ of $t'$
that contains the origin.  We then grow outwards, so that the tiling
$t'$ is combinatorially identical to $t$, only with each tile of type
$C_j$ replaced with a tile of type $C'_j$, and each edge of type $I_j$
replaced by $I'_j$.  This is shown in Figure 2, where a patch of the
original Penrose tiling is replaced by a patch of rational Penrose
tiling.

\begin{figure}
\vbox{
%\centerline{\epsfysize=9.5cm \epsfbox{ssft/R180} }
\centerline{\epsfysize=9.5cm \epsfbox{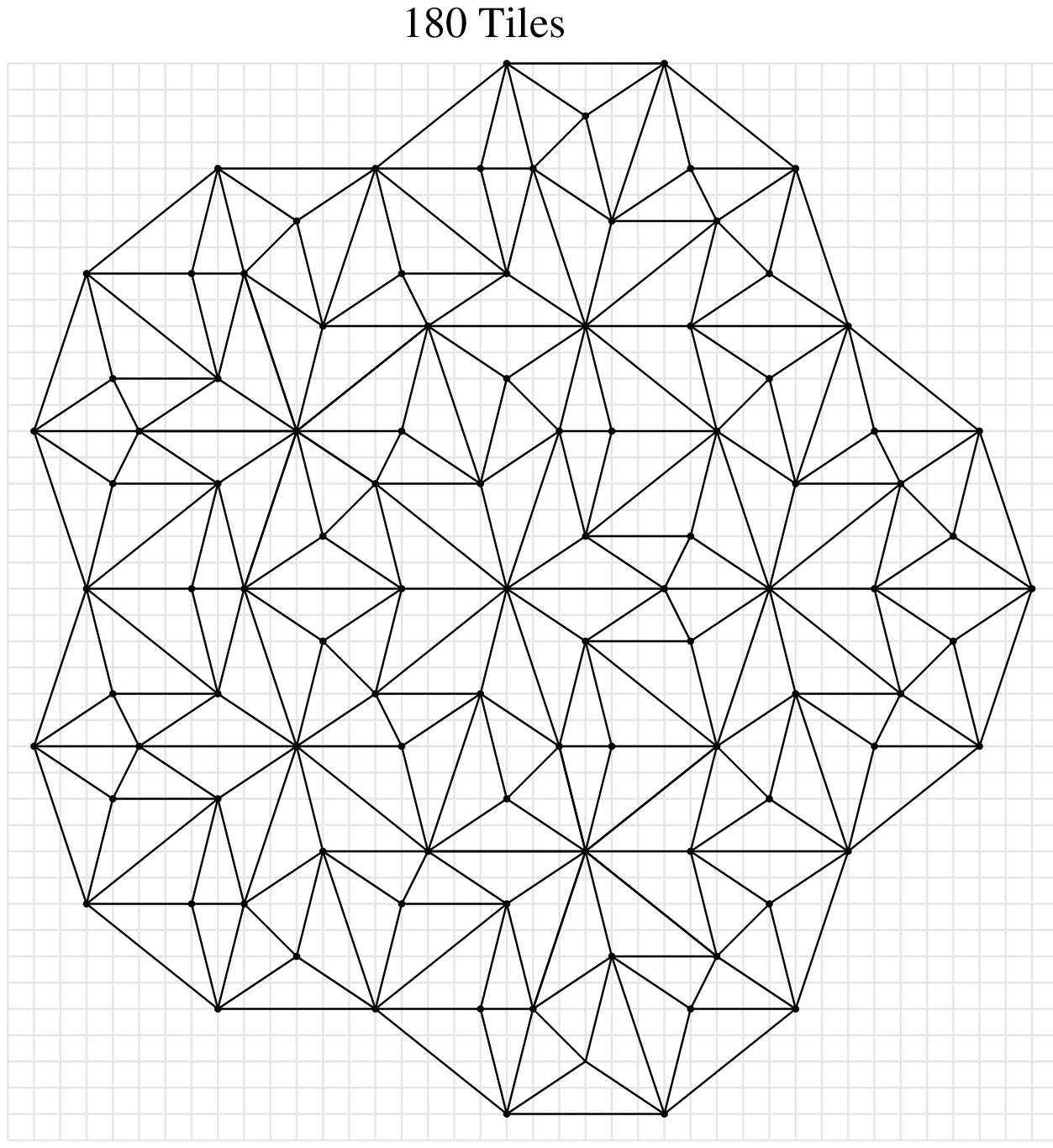} }
}
\caption{A patch of the rational Penrose tiling}
\end{figure}

To see that this construction does result in a tiling, we must show
that the vertices of $t'$ are well defined.  Let $x$ be a vertex of
$t$, and consider two paths from a vertex $y$ of the central seed tile
to $x$.  The algebraic difference of these two paths, namely zero, is
the boundary of a sum of tiles in $t$.  By equations (\ref{eq5}), the
algebraic difference of the corresponding sum of vectors $v'$ is also
zero.  This means that either path can be used to determine the
position of $x'$, the vertex in $t'$ that corresponds to $x$.  Once
the vertices are defined, the edges and faces follow.

This transformation is continuous.  If two tilings $t$ and $\tilde t$
agree on a large neighborhood of the origin, then $t'$ and $\tilde t'$
agree on a large neighborhood of the origin.  If $t$ and $\tilde t$
differ by a small translation, then $t'$ and $\tilde t'$ differ by a 
small translation, as determined by the homeomorphism between the 
center tile of $t$ and that of $t'$. (As noted above, there is no 
ambiguity, and no discontinuity, if the origin in $t$ 
sits on the boundary of a tile.) Similarly, the reverse transformation,
from tilings in $R$ to tilings in $P$, is also continuous.  Thus $P$ 
and $R$ are homeomorphic tiling spaces.

Finally, since each tiling $t$ in $P$ is combinatorially equivalent to
a tiling $t'$ in $R$, any local atlas for the $P$ system can be
naturally transformed into a local atlas for the $R$ system, and
vice-versa.

In higher dimensions, the argument is essentially as in two dimensions.
As mentioned above, since we are taking our tiles in higher dimensions
to be {\it triangulated }, we can proceed with a finite system of
homogeneous equations---all of the form $a+b+c=0$. As before, we can 
find a rational solution arbitrarily close to the original 
vectors.  To construct homeomorphisms between tiles $C_k$ and $C_k'$, 
one must start with homeomorphisms (e.g. linear maps) between edges 
$I_i$ and $I_i'$, extend these to homeomorphisms of the 2-skeleton, then
of the 3-skeleton, and so on.  There are no topological obstructions.  
\QED

\enddemo

To complete the proof of Theorem 1, we must only prove

\begin{lem} Every rational tiling space is a fiber bundle over the torus.
\end{lem}

\demo{Proof.} Let $R$ be a rational tiling, and let $D$ be the least 
common multiple of all the denominators
of all the coordinates of displacement vectors $v_i$ for the tiles in $R$.
Rescale $R$ by $D$, so that all displacement vectors are integers.  Then
all the vertices in any fixed tiling have the same coordinates (mod $\zed^d$).
These coordinates give a natural projection from the space of tilings 
to the $d$-torus $\red^d/\zed^d$. \QED
\enddemo
 
\section{Square tiling spaces}

We have shown that our general tiling space $P$ is homeomorphic to a
rational tiling space $R$ that is of finite type if $P$ is (and is not
if $P$ is not).  By rescaling, we can assume that $R$ is
in fact integral.  Topological conjugacies preserve finite type
[RS].  To complete the proof of Theorem 2, it suffices to prove

\begin{lem} Every rational tiling space $R$ is, after rescaling, 
topologically conjugate to a square-type tiling space $S$. 
\end{lem}

\demo{Proof.} As before, we work first in 2 dimensions, and then
sketch what modifications need to be made in higher dimensions.  Also
as before, we illustrate each step of our construction with the
Penrose system.

We first rescale $R$ so that $R$ becomes an integral tiling.
Furthermore, we assume that each tile contains a circle of radius
greater than $\sqrt{2}/2$; this can always be achieved by further
scaling.  Next we replace each of our straight edges $I'_i$ with
zig-zags $J_i$, that is with sequences of unit displacements in the
coordinate directions.  We do this in such a way that the maximum
distance of a point in $J_i$ from the original edge $I'_i$ is
minimized.  In particular, one can always choose $J_i$ such that this
distance is no greater than $\sqrt{2}/2$. There is sometimes more than
one way to minimize this distance.  For example, one could replace a
diagonal edge from (0,0) to (1,1) with a zig-zag from (0,0) to (1,0)
to (1,1), or with a zig-zag from (0,0) to (1,0) to (1,1).  In such a case,
one must make a choice and apply it consistently.  A possible
set of zig-zags for the rational Penrose system is given in Figure
3. Under this replacement, the 180-tile patch of Figure 2 turns into
the patch of Figure 4.

\begin{figure}
\vbox{
%\centerline{\epsfysize=9.5cm \epsfbox{ssft/edges} }}
\centerline{\epsfysize=9.5cm \epsfbox{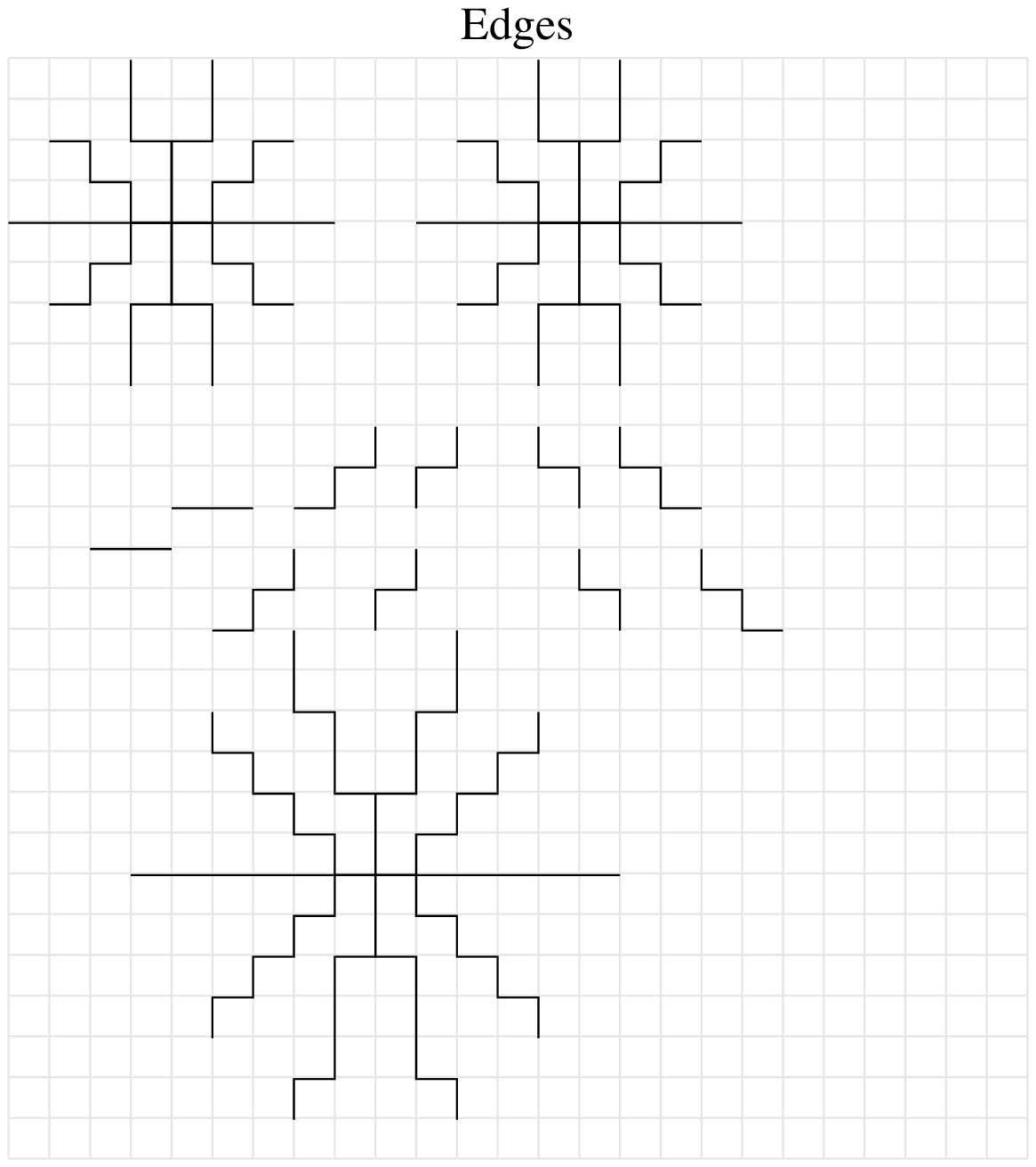} }}
\caption{The 40 zig-zag edges.  The ``a'' edges are in the upper left, the
``b'' edges in the upper right, the ``c'' edges in the middle, and the ``d''
edges at the bottom.}
\end{figure}

\begin{figure}
\vbox{
%\centerline{\epsfysize=9.5cm \epsfbox{ssft/S180} }}
\centerline{\epsfysize=9.5cm \epsfbox{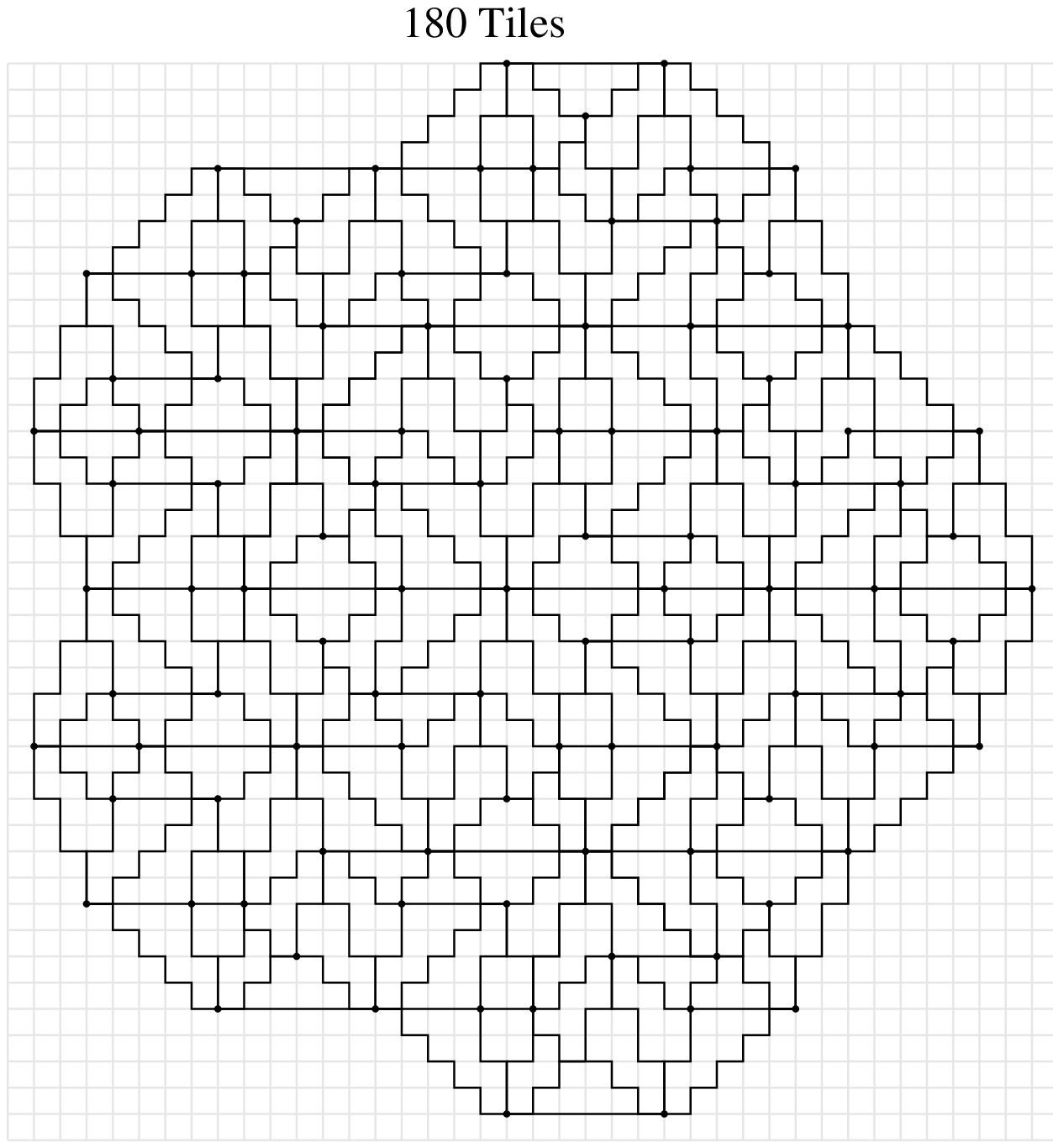} }}
\caption{A patch of the square Penrose tiling}
\end{figure}

This defines a space $Z$ of tilings whose edges are zig-zags. To each
tiling $t'$ in $R$ we generate a tiling $z$ in $Z$ by replacing each
straight edge in $t'$ with its corresponding zig-zag.  If a tile type
$C'_k$ in the $R$ system is bounded by several straight edges $I'_i$,
then the tile type $D_k$ in the $Z$ system is defined to be the region
bounded by the corresponding zig-zags $J_i$'s. The condition that
$C'_k$ contains a circle of radius greater that $\sqrt{2}/2$ ensures
that $D_k$ is nonempty.  (It may, however, be disconnected).  It may
happen that geometrically non-congruent tile types $C'_k$ generate congruent
tile types $D_k$; however, as marked tiles, these $D_k$'s should be 
considered distinct.

The operation of replacing straight edges with zig-zags is reversible
and does not require a choice of origin.  It therefore commutes with
translation and defines a topological conjugacy between $R$ and $Z$. 

In the tiling system $Z$, the basic tiles are irregularly shaped regions 
$D_k$ bounded by zig-zags, and we have already seen that each $D_k$ is
nonempty.  
Suppose that the tile $D_k$ has area $n$. Then $D_k$ can be decomposed as
the union of $n$ unit squares $D_k^1, \ldots, D_k^n$.  In the tiling system
$S$, the basic tiles are the squares $D_k^i$, and we apply a matching rule
that says that wherever one of the $D_k^i$ squares is found, 
the other $n-1$ squares that make up $D_k$ are also found nearby, arranged
to form the larger region $D_k$.  A tiling in $S$ can therefore be amalgamated
into a tiling by tiles $D_k$.  We allow in $S$ those tilings, and only
those tilings, that amalgamate into tilings in $Z$.  In this way, the 
tiling system $S$ is naturally conjugate to $Z$. 

The proof in higher dimensions is almost identical.
In dimension 3, one must pick zig-zags $J_i$ to replace
the straight edges $I'_i$.  If several edges $I'_i$ bound a 2-face of a 
tile $C'_k$, we must find a union of unit squares (oriented in the coordinate
directions), bounded by the appropriate zig-zags $J_i$, that approximates
this face.  The tile $D_k$ is then the solid region bounded by these 
zig-zag faces.  In dimension $d>3$, one works recursively, replacing
edges $I_i'$ with zig-zags $J_i$, then replacing 2-cells with unions of 
squares, 3-cells with unions of cubes, and so on up through dimension $d-1$.
The tiles $D_k$ are the $d$-cells bounded by the $d-1$ cells constructed
in this manner.  One can compute a universal bound for each dimension, so
that the $d-1$ dimensional zig-zags are within that universal bound of the
original faces of the $C'_k$'s.  As long as the $C'_k$'s contain a sphere
of radius greater than that bound, the resulting $D_k$ will be nonempty. 
\QED

\enddemo

\bigskip

\noindent
\author{
Lorenzo Sadun and R.F. Williams\\
Department of Mathematics\\
The University of Texas at Austin\\
Austin, TX 78712 U.S.A.\\
sadun@math.utexas.edu and bob@math.utexas.edu}


\begin{thebibliography}{HW}

\bibitem[AP]{AP} Jared Anderson and Ian Putnam,
Topological invariants for substitution tilings and their associated C*
algebras, 
Ergodic Theory and Dynamical Systems {\bf 18} part 3 (1998) 509--538.

\bibitem[B]{B} M. Boyle, Symbolic dynamics and matrices, IMA Vol. in
Math. and Its Appl. {\bf 50}, Springer-Verlag, 1-38.

\bibitem[BF]{BF} R. Bowen and J. Franks, Homology for zero dimensional
basic sets, Annals of Math., {\bf 106}(1977), 73-92. 


\bibitem[BJV]{BJV} Marcy Barge, James Jacklitch, and Gioia Vago,
Homeomorphisms of one-dimensional inverse limits with application to
substitution tilings, unstable manifolds, and tent maps,
Contemporary Mathematics, {\bf 246}, 1 - 15.

\bibitem[F]{F} J. Franks, Flow equivalence of subshifts of finite type,
Ergodic Theory and Dynamical Systems {\bf 4} (1984) 53-66. 

\bibitem[FJ]{FJ} F. T. Farrell and L. E. Jones,  
New attractors in hyperbolic dynamics, 
Journal of Differential Geometry {\bf 15}, no. 1 (1980), 107--133.

\bibitem[LM]{LM} D. Lind and B. Marcus, 
``An Introduction to Symbolic Dynamics and Coding'', Cambridge  Univ.
Press, Cambridge 1995

\bibitem[ORS]{ORS} Nicolus Ormes, Charles Radin, and Lorenzo Sadun, 
A homeomorphism invariant for substitution tiling spaces,
preprint, University of Texas at Austin.

\bibitem[PS]{PS} B. Parry and D. Sullivan, A topological invariant for
  flows on 1-dimensional spaces, Topology {\bf 14},(1975) 297-299.

\bibitem[R]{R} C. Radin, The pinwheel tilings of the plane, Annals of
Math.  {\bf 139} (1994), 661--702.

\bibitem[RS]{RS} C. Radin and L. Sadun, 
Isomorphisms of Hierarchical Structures, to appear in
Ergodic Theory \& Dynamical Systems. 

\bibitem[S]{S} M. Sullivan, Invariants of twist-wise flow equivalence,
  Discrete and Continuous Dynamical Systems {\bf 4}(1998) 474-484.

\bibitem[Wa]{Wa} Proceedings on Warwick Symposium on $\zed^n$ Actions, 
LMS Lecture Notes Series {\bf 228} Cambridge Univ. Press, 1996.

\bibitem[W1]{W1} 
R. F. Williams
One dimensional nonwandering sets, 
Topology {\bf 6},(1967), 473--487.

\bibitem[W2]{W2} 
R. F. Williams,
Expanding attractors, 
Institute des Hautes \'Etudes
Scientifique Publ. Math. no. 43(1973), 169--203.

\bibitem[W3]{W3} R. F. Williams,
The Penrose, Ammann and DA tiling spaces are Cantor set fiber bundles.
To appear in Ergodic theory and dynamical systems.
See the home page of R. F. Williams.
\end{thebibliography}
\end{document}